\documentclass[12pt,a4paper]{amsart}
\usepackage{mathrsfs}

\newtheorem{theo+}              {Theorem}           [section]
\newtheorem{prop+}  [theo+]     {Proposition}
\newtheorem{coro+}  [theo+]     {Corollary}
\newtheorem{lemm+}  [theo+]     {Lemma}
\newtheorem{exam+}  [theo+]     {Example}
\newtheorem{rema+}  [theo+]     {Remark}
\newtheorem{defi+}  [theo+]     {Definition}

\newenvironment{theorem}{\begin{theo+}}{\end{theo+}}
\newenvironment{proposition}{\begin{prop+}}{\end{prop+}}
\newenvironment{corollary}{\begin{coro+}}{\end{coro+}}
\newenvironment{lemma}{\begin{lemm+}}{\end{lemm+}}
\newenvironment{definition}{\begin{defi+}}{\end{defi+}}

\usepackage{amsthm}
\theoremstyle{plain} \theoremstyle{remark}

\renewcommand{\Bbb}{\mathbb}
\newcommand{\rank}{\mbox{rank}}

\def\E{/\kern-1.0em \equiv }
\def\h{harmonic morphism}
\def\bih{biharmonic morphism}

\title{ Biharmonic maps into Sol and Nil spaces  }
\author{Ye-Lin Ou and Ze-Ping Wang}

\address{Department of Mathematics,\newline\indent Texas A $\&$ M University-Commerce,
\newline\indent Commerce TX 75429, USA.\newline\indent
E-mail:yelin$\_$ou@tamu-commerce.edu\;(Ou)\newline\indent\newline\indent
School of Math and Computer Science,\newline\indent Guangxi
University for Nationalities,\newline\indent Nanning 530006, P. R.
China. }

\begin{document}

\title[ Biharmonic maps into Sol and Nil spaces]{ Biharmonic maps into Sol and Nil spaces  }
 \maketitle

\section*{Abstract}
\begin{quote}
{\footnotesize In this paper, we study biharmonic maps into Sol and
Nil spaces, two model spaces of Thurston's 3-dimensional geometries.
We characterize non-geodesic biharmonic curves in Sol space and
prove that there exists no non-geodesic biharmonic helix in Sol
space. We also show that a linear map from a Euclidean space into
Sol or Nil space is biharmonic if and only if it is a harmonic map,
and give a complete classification of such maps.}
\end{quote}
\section{introduction}

In this paper, we work in the category of smooth objects, so
manifolds, maps, vector fields, etc, are assumed to be smooth unless it is stated otherwise.\\

 A map $\varphi:(M, g)\longrightarrow (N, h)$ between
Riemannian manifolds is called a {\em biharmonic map} if
$\varphi|\Omega$ is a critical point of the bienergy
\begin{equation}\nonumber
E^{2}\left(\varphi,\Omega \right)= \frac{1}{2} {\int}_{\Omega}
\left|\tau_{2}(\varphi) \right|^{2}{\rm d}x
\end{equation}
for every compact subset $\Omega$ of $M$, where
$\tau_{2}(\varphi)={\rm Trace}_{g}\nabla {\rm d} \varphi$ is the
tension field of $\varphi$. Using the first variational formula (see
\cite{Ji}) one sees that $\varphi$ is a biharmonic map if and only
if its bitension field vanishes identically, i.e.,
\begin{equation}\label{BI1}
\tau^{2}(\varphi):=-\triangle^{\varphi}(\tau_{2}(\varphi)) - {\rm
Trace}_{g} R^{N}({\rm d}\varphi, \tau_{2}(\varphi)){\rm d}\varphi
=0,
\end{equation}
where
\begin{equation}\notag
\triangle^{\varphi}=-{\rm Trace}_{g}(\nabla^{\varphi})^{2}= -{\rm
Trace}_{g}(\nabla^{\varphi}\nabla^{\varphi}-\nabla^{\varphi}_{\nabla^{M}})
\end{equation}
is the Laplacian on sections of the pull-back bundle $\varphi^{-1}
TN$ and $R^{N}$ is the curvature operator of $(N, h)$ defined by
$$R^{N}(X,Y)Z=
[\nabla^{N}_{X},\nabla^{N}_{Y}]Z-\nabla^{N}_{[X,Y]}Z.$$ Note that
$\tau^{2}(\varphi)=-J^{\varphi}(\tau_{2}(\varphi))$, where
$J^{\varphi}$ is the Jacobi operator which plays an important role
in the study of harmonic maps.\\
\indent Clearly, any harmonic map is  biharmonic, so harmonic maps
are a subclass of biharmonic maps. It is interesting to note that so
far, apart from maps between Euclidean spaces (in which case any
polynomial map of degree less than four is a biharmonic map), not
many example of {\em proper} (meaning non-harmonic) biharmonic maps
between Riemannain manifolds have been found (see, e.g., \cite{MO},
\cite{LOU}, \cite{Ou1}, and the bibliography of biharmonic maps
\cite{LMO}).
\\

Recently, some work has been done in the study of non-geodesic
biharmonic curves in some model spaces. For example, \cite{CMO1}
gives a complete classification of non-geodesic biharmonic curves in
$3$-sphere while in \cite{COP} it is proved that a non-geodesic
biharmonic curve in Heisenberg group $\mathbb{H}_{3}$ is a helix and
the explicit parametrizations of such curves are given. For the
study of biharmonic curves in Berger's spheres, in Minkowski
3-space, in Cartan-Vranceanu 3-dimensional space, and in contact and
Sasakian manifolds see \cite{Ba},
\cite{In}, \cite{CMOP}, and \cite{EY} respectively.\\

In this paper, we first characterize non-geodesic biharmonic curves
in Sol space (Theorem \ref{PE010}) and show that there exists no
non-geodesic biharmonic helix in Sol space (Theorem \ref{PE0105}).
In the second part of the paper, we study linear biharmonic maps
from Euclidean space $\mathbb{R}^{m}$ into Sol space
$(\mathbb{R}^{3}, g_{Sol})$ and Nil space $(\mathbb{R}^{3},
g_{Nil})$ using the linear structure of the underlying manifolds of
Sol and Nil spaces. We show that such a linear map is biharmonic if
and only if it is a harmonic
map and give a complete description of such maps (Theorems \ref{MAI1} and \ref{MAI2}).\\

\section{Biharmonic curves in Sol space}

{\bf Biharmonic curve equation in Frenet Frames}. Let
$I\subset\mathbb{R}$ be an open interval, and $\gamma:
I\longrightarrow(N,h)$ be a curve, parametrized by arc length, on a
Riemannian manifold. Putting $T=\gamma'$, we can write the tension
field of $\gamma$ as $\tau(\gamma)=\nabla_{\gamma'}{\gamma'}$, and
the biharmonic map equation (\ref{BI1}) reduces to
\begin{equation}\notag
 \nabla^{3}_{T}T-R^{N}(T,\nabla_{T}T)T=0.
\end{equation}

A successful key to study the geometry of a curve is to use the
Frenet frames along the curve which is recalled in the following
\begin{definition}\label{PE0101}(see, for example, \cite{La})\\
The Frenet frame $\{F_{i}\}_{i=1, 2, \ldots, n}$ associated to a
curve $\gamma: I\subset \mathbb{R}\longrightarrow(N^{n}, h)$
parametrized by arc length is the orthonormalisation of the
(n+1)-tuple\\ $\{\nabla_{\frac{\partial}{\partial
t}}^{(k)}d\gamma(\frac{\partial}{\partial t})\}_{k=0, 1,\; 2,\;
\ldots,\; n}$ described by:

\begin{equation}\notag
\begin{array}{lll}
F_{1}=d\gamma(\frac{\partial}{\partial t})\\
\nabla_{\frac{\partial}{\partial t}}^{(\gamma)}F_{1}=k_{1} F_{2}\\
\nabla_{\frac{\partial}{\partial t}}^{(\gamma)}F_{i}=-k_{i-1}
F_{i-1}+k_{i} F_{i+1},\forall\; i=2,\;3,\;\ldots,\;n-1\\
\nabla_{\frac{\partial}{\partial t}}^{(\gamma)}F_{n}=-k_{n-1}
F_{n-1},\\
\end{array}
\end{equation}
where the functions $\{k_{1}, k_{2}, \ldots, k_{n-1}\}$ are called
the curvatures of $\gamma$ and $\nabla^{\gamma}$ is the connection
on the pull-back bundle $\gamma^{-1}(TN)$. Note that
$F_{1}=T=\gamma'$ is the unit tangent vector field along the curve.
\end{definition}
With respect to the Frenet frame, the biharmonic curve equation
takes the following form.
\begin{lemma}\label{PE0102}(see, for example, \cite{CMOP})
Let $\gamma: I\subset \mathbb{R}\longrightarrow(N^{n}, h) (n\geq 2)$
be a curve parametrized by arc length from an open interval of
$\mathbb{R}$ into a Riemannian manifold $(N^{n},h)$. Then $\gamma$
is biharmonic if and only if :
\begin{equation}\label{N3301}
\left\{\begin{array}{rl}
 k_{1}k_{1}'=0\\
k_{1}''-k_{1}^{3}-k_{1}k_{2}^{2}+k_{1}R(F_{1},F_{2},F_{1},F_{2})=0\\
2k_{1}'k_{2}+k_{1}k_{2}'+k_{1}R(F_{1},F_{2},F_{1},F_{3})=0\\
k_{1}k_{2}k_{3}+k_{1}R(F_{1},F_{2},F_{1},F_{4})=0\\
k_{1}R(F_{1},F_{2},F_{1},F_{j})=0\;\;\;\;\;\;j=5,\ldots,n.\;\;
\end{array}\right.
\end{equation}
\end{lemma}

 {\bf Riemannian structure of  Sol space}. Sol space, one of Thurston's eight
 3-dimensional geometries, can be viewed as $\mathbb{R}^{3}$
 provided with Riemannian metric $g_{Sol}={\rm d}s^{2}=e^{2z}{\rm d}x^{2}+e^{-2z}{\rm d}y^{2}+{\rm
 d}z^{2}$, where $(x, y, z)$ are the standard coordinates in
 $\mathbb{R}^{3}$. Note that the Sol metric can also be written as:
\begin{align}\nonumber
ds^{2}=\sum_{i=1}^{3}\omega^{i}\otimes\omega^{i},
\end{align}
where
\begin{equation}\notag
\omega^{1}=e^{z}dx,\;\omega^{2}=e^{-z}dy,\;\omega^{3}=dz,
\end{equation}
and the orthonormal basis dual to the 1-forms is
\begin{equation}\notag
e_{1}=e^{-z}\frac{\partial}{\partial x},\;
e_{2}=e^{z}\frac{\partial}{\partial y},
\;e_{3}=\frac{\partial}{\partial z}.
\end{equation}
With respect to this orthonormal basis, the Levi-Civita connection
and the Lie brackets can be easily computed as:
\begin{equation}\label{00E220}
\begin{array}{lll}
\nabla_{e_{1}}e_{1}=-e_{3},\nabla_{e_{1}}e_{2}=0,\nabla_{e_{1}}e_{3}=e_{1}\\
\nabla_{e_{2}}e_{1}=0,\nabla_{e_{2}}e_{2}=e_{3},\nabla_{e_{2}}e_{3}=-e_{2}\\
\nabla_{e_{3}}e_{1}=0,\nabla_{e_{3}}e_{2}=0,\nabla_{e_{3}}e_{3}=0.\\
\end{array}
\end{equation}
\begin{equation}\label{00E23}
[e_{1},e_{2}]=0, \;[e_{2},e_{3}]=-e_{2}, \;[e_{1},e_{3}]=e_{1}.
\end{equation}
We adopt the following notation and sign convention for Riemannian
curvature operator.
\begin{equation}\label{R1}
 R(X,Y)Z=\nabla_{X}\nabla_{Y}Z
-\nabla_{Y}\nabla_{X}Z-\nabla_{[X,Y]}Z,\\
\end{equation}
 the Riemannian curvature tensor is
given by
\begin{equation}\label{R2}
\begin{array}{lll}
 R(X,Y,Z,W)=g( R(Y,X)Z,W)=-g( R(X,Y)Z,W).\\

\end{array}
\end{equation}
 Moreover we
put
\begin{equation}\label{R3}
R_{abc}=R(e_{a},e_{b})e_{c}, \; R_{abcd}=R(e_{a},
e_{b},e_{c},e_{d}),
\end{equation}
where  the indices $ a, b, c, d$ take the values $1, 2, 3$. \\

A direct computation using (\ref{00E220}), (\ref{00E23}),
(\ref{R1}), (\ref{R2}), and (\ref{R3}) gives the following non-zero
components of Riemannian curvature of Sol space with respect to the
orthonomal basis $\{ e_{1}, e_{2}, e_{3}\}$ (we do not list those
that can be obtained by symmetric properties of curvature):
\begin{equation}\notag
\begin{array}{lll} R_{121}=-e_{2},
R_{131}=e_{3}, R_{122}=e_{1},\\R_{232}=e_{3},
R_{133}=-e_{1},R_{233}=-e_{2},\\
\end{array}
\end{equation}
and
\begin{equation}\label{00E24}
\begin{array}{lll}
 R_{1212}=-g(R(e_{1},e_{2})e_{1},e_{2})=-g(-e_{2},e_{2})=1,\\
R_{1313}=-g(R(e_{1},e_{3})e_{1},e_{3})=-g(e_{3},e_{3})=-1,\\
R_{2323}=-g(R(e_{2},e_{3})e_{2},e_{3})=-g(e_{3},e_{3})=-1.
\end{array}
\end{equation}

{\bf Biharmonic curves in Sol space}. Let $\gamma: I\longrightarrow
(\mathbb{R}^{3}, g_{Sol})$ be a curve on Sol space parametrized by
arc length. Let $\{T, N, B\}$ be the Frenet frame fields tangent to
Sol space along $\gamma$ defined as follows: $T$ is the unit vector
field $\gamma'$ tangent to $\gamma$, $N$ is the unit vector field in
the direction of $\nabla_{T}T $ (normal to $\gamma$), and $B$ is
chosen so that $\{T, N, B\}$ is a positively oriented orthonormal
basis. Then, by Definition (\ref{PE0101}), we have the following
Frenet formulas

\begin{equation}\label{00E26}
\left\{\begin{array}{rl}
 \nabla_{T}T=kN\\
 \nabla_{T}N=-kT+\tau B\\
 \nabla_{T}B=-\tau N,
\end{array}\right.
\end{equation}
where $k=|\nabla_{T}T | $ is the geodesic curvature
 and $\tau $ the geodesic torsion of $\gamma$. With respect to the orthonormal basis
 $\{ e_{1}, e_{2}, e_{3}\}$ we can write $T= T_{1}e_{1}+T_{2}e_{2}+T_{3}e_{3}, \;N=
N_{1}e_{1}+N_{2}e_{2}+N_{3}e_{3}$, $B=T\times N=
B_{1}e_{1}+B_{2}e_{2}+B_{3}e_{3}$, and we have

\begin{theorem}\label{PE010}
Let $\gamma: I\longrightarrow (\mathbb{R}^{3}, g_{Sol})$ be a curve
parametrized by arc length. Then $\gamma$ is a non-geodesic
biharmonic curve if and only if
\begin{equation}\label{00E27}
\left\{\begin{array}{rl}

 k={\rm constant}\neq 0\\
k^{2}+\tau^{2}=2B_{3}^{2}-1\\
\tau'=2N_{3}B_{3}.\\
\end{array}\right.
\end{equation}
\end{theorem}
\begin{proof}
By (\ref {N3301}) of lemma (\ref {PE0102}) we see that $\gamma$ is a
biharmonic curve if and only if
\begin{equation}\notag
\left\{\begin{array}{rl}
 kk'=0\\
k''-k^{3}-k\tau^{2}+kR(T,N,T,N)=0\\
2k'\tau+k\tau'+kR(T,N,T,B)=0,\\
\end{array}\right.
\end{equation}
which is equivalent to
\begin{equation}\label{N3303}
\left\{\begin{array}{rl}
 k=constant\neq 0\\
k^{2}+\tau^{2}=R(T,N,T,N)\\
\tau'=-R(T,N,T,B)\\
\end{array}\right.
\end{equation}
since $k\neq 0$ by the assumption that $\gamma$ is non-geodesic. A
direct computation using (\ref{00E24}) yields
\begin{eqnarray*}
\begin{array}{lll}
R(T,N,T,N)= \sum_{i,j,l,p=1}^{3}T_{l}N_{p}T_{i}N_{j}R_{lpij}\\
=T_{1}N_{2}T_{1}N_{2}R_{1212}+T_{1}N_{2}T_{2}N_{1}R_{1221}
+T_{2}N_{1}T_{2}N_{1}R_{2121}+T_{2}N_{1}T_{1}N_{2}R_{2112}\\
+T_{1}N_{3}T_{1}N_{3}R_{1313}+T_{1}N_{3}T_{3}N_{1}R_{1331}
+T_{3}N_{1}T_{3}N_{1}R_{3131}+T_{3}N_{1}T_{1}N_{3}R_{3113}\\
+T_{2}N_{3}T_{2}N_{3}R_{2323}+T_{2}N_{3}T_{3}N_{2}R_{2332}
+T_{3}N_{2}T_{3}N_{2}R_{3232}+T_{3}N_{2}T_{2}N_{3}R_{3223}\\
=T_{1}^{2}N_{2}^{2}-T_{1}T_{2}N_{1}N_{2}+T_{2}^{2}N_{1}^{2}-T_{1}T_{2}N_{1}N_{2}\\
-T_{1}^{2}N_{3}^{2}+T_{1}T_{3}N_{1}N_{3}-T_{3}^{2}N_{1}^{2}+T_{1}T_{3}N_{1}N_{3}\\
-T_{2}^{2}N_{3}^{2}+T_{2}T_{3}N_{2}N_{3}-T_{3}^{2}N_{2}^{2}+T_{2}T_{3}N_{2}N_{3}\\
=(T_{1}N_{2}-T_{2}N_{1})^{2}-(T_{3}N_{1}-T_{1}N_{3})^{2}-(T_{3}N_{2}-T_{2}N_{3})^{2}\\
=B_{3}^{2}-B_{1}^{2}-B_{2}^{2} =2B_{3}^{2}-1\\\notag (by\; T\times
N=B ,T\times B=-N,B_{1}^{2}+B_{2}^{2}+B_{3}^{2}=1),
\end{array}
\end{eqnarray*}
and
\begin{eqnarray*}
\begin{array}{lll}
R(T,N,T,B)=\sum_{i,j,l,p=1}^{3}T_{l}N_{p}T_{i}B_{j}R_{lpij}\\
=T_{1}N_{2}T_{1}B_{2}R_{1212}+T_{1}N_{2}T_{2}B_{1}R_{1221}\\
+T_{2}N_{1}T_{2}B_{1}R_{2121}+T_{2}N_{1}T_{1}B_{2}R_{2112}\\\notag
+T_{1}N_{3}T_{1}B_{3}R_{1313}+T_{1}N_{3}T_{3}B_{1}R_{1331}\\
+T_{3}N_{1}T_{3}B_{1}R_{3131}+T_{3}N_{1}T_{1}B_{3}R_{3113}\\\notag
+T_{2}N_{3}T_{2}B_{3}R_{2323}+T_{2}N_{3}T_{3}B_{2}R_{2332}\\
+T_{3}N_{2}T_{3}B_{2}R_{3232}+T_{3}N_{2}T_{2}B_{3}R_{3223}\\\notag
=T_{1}^{2}N_{2}B_{2}-T_{1}T_{2}N_{1}B_{2}+T_{2}^{2}N_{1}B_{1}-T_{1}T_{2}N_{2}B_{1}\\\notag
-T_{1}^{2}N_{3}B_{3}+T_{1}T_{3}N_{1}B_{3}-T_{3}^{2}N_{1}B_{1}+T_{1}T_{3}B_{1}N_{3}\\\notag
-T_{2}^{2}N_{3}B_{3}+T_{2}T_{3}N_{2}B_{3}-T_{3}^{2}N_{2}B_{2}+T_{2}T_{3}B_{2}N_{3}\\\notag
=(T_{1}B_{2}-T_{2}B_{1})(T_{1}N_{2}-T_{2}N_{1})\\\notag
-(T_{3}N_{1}-T_{1}N_{3})(T_{3}B_{1}-T_{1}B_{3})\\\notag
-(T_{3}N_{2}-T_{2}N_{3})(T_{3}B_{2}-T_{2}B_{3})\\\notag
=-N_{3}B_{3}+N_{2}B_{2}+N_{1}B_{1} =-2N_{3}B_{3}\\\notag (by\;
T\times N=B ,T\times B=-N, N_{1}B_{1}+N_{1}B_{2}+N_{1}B_{3}=0),
\end{array}
\end{eqnarray*}
these, together with Equation (\ref{N3303}), complete the proof of
the theorem.
\end{proof}

As an immediate consequence we have
\begin{corollary}\label{PE0}
Let $\gamma:I\longrightarrow (\mathbb{R}^{3}, g_{Sol})$ be a
non-geodesic curve parametrized by arc length. If $B_{3}= 0$, then
$\gamma$ is not biharmonic.
\end{corollary}
\begin{corollary}\label{PE0103}
Let $\gamma: I\longrightarrow (\mathbb{R}^{3}, g_{Sol})$ be a
non-geodesic curve parametrized by arc length. If $B_{3}$ is
constant and  $N_{3}B_{3}\neq 0$, then $\gamma$ is not biharmonic.
\end{corollary}

Similar to the terminology used for curves in $\mathbb{R}^{3}$, we
keep the name helix for a curve in a Riemannian 3-manifold having
constant both geodesic curvature and geodesic torsion. With this
terminology, we can use Equation (\ref{00E27}) to deduce the
following
\begin{corollary}\label{PE0104}
Let $\gamma:I\longrightarrow (\mathbb{R}^{3}, g_{Sol})$ be a
non-geodesic biharmonic helix parametrized by arc length, then
\begin{equation}\label{00E32}
\left\{\begin{array}{rl}

B_{3}={\rm constant}\neq 0\\
N_{3}=0\\
k^{2}+\tau^{2}=2B_{3}^{2}-1.\\
\end{array}\right.
\end{equation}
\end{corollary}

Non-geodesic biharmonic helices in 3-dimensional sphere $S^3$, in
Heisenberg group $H_3$, and in Cartan-Vranceanu 3-dimensional space
have been studied in \cite{CMO1}, \cite{COP}, and \cite{CMOP}
respectively. In contrast to the situations in those 3-dimensional
spaces where there are rich examples of such curves our next theorem
shows that there exists no such curves in Sol space.

\begin{theorem}\label{PE0105}
There exists no non-geodesic biharmonic helix in Sol space.
\end{theorem}

\begin{proof}
Suppose that $\gamma:I\longrightarrow (\mathbb{R}^{3}, g_{Sol})$ is
a non-geodesic biharmonic helix parametrized by arc length. We shall
derive a contradiction by showing that $\gamma$ must be a geodesic.
We can use (\ref{00E220}) to compute the covariant derivatives  of
the vector fields $T, N, B$ as:
\begin{equation}\label{00E34}
\left\{\begin{array}{rl}
\nabla_{T}T=(T_{1}'+T_{1}T_{3})e_{1}+(T_{2}'-T_{2}T_{3})e_{2}+(T_{2}^{2}-T_{1}^{2}+T_{3}')e_{3}\\
\nabla_{T}N=(N_{1}'+T_{1}N_{3})e_{1}+(N_{2}'-T_{2}N_{3})e_{2}+(T_{2}N_{2}-T_{1}N_{1}+N_{3}')e_{3}\\
\nabla_{T}B=(B_{1}'+T_{1}B_{3})e_{1}+(B_{2}'-T_{2}B_{3})e_{2}+(T_{2}B_{2}-T_{1}B_{1}+B_{3}')e_{3}.\\
\end{array}\right.
\end{equation}
It follows that the third components of these vectors are given by
\begin{equation}\label{00E35}
\left\{\begin{array}{rl}
\langle\nabla_{T}T,e_{3}\rangle=T_{2}^{2}-T_{1}^{2}+T_{3}'\\
\langle\nabla_{T}N,e_{3}\rangle=T_{2}N_{2}-T_{1}N_{1}+N_{3}'\\
\langle\nabla_{T}B,e_{3}\rangle=T_{2}B_{2}-T_{1}B_{1}+B_{3}'.\\
\end{array}\right.
\end{equation}

On the other hand, using Frenet formulas (\ref {00E26}) we have,
\begin{equation}\label{00AE35}
\left\{\begin{array}{rl}
\langle\nabla_{T}T,e_{3}\rangle=kN_{3}\\
\langle\nabla_{T}N,e_{3}\rangle=-kT_{3}+\tau
B_{3}\\
\langle\nabla_{T}B,e_{3}\rangle=-\tau
N_{3}.\\
\end{array}\right.
\end{equation}
Since $\gamma$ is assumed to be a non-geodesic biharmonic helix, we
have, by Corollary \ref{PE0104},
 $N_3=0,\; B_3={\rm constant}$. These,
 together with Equations (\ref{00E35}) and (\ref{00AE35}), give

\begin{equation}\label{00E36}
\left\{\begin{array}{rl}
T_{2}^{2}-T_{1}^{2}+T_{3}'=0\\
T_{2}N_{2}-T_{1}N_{1}=-kT_{3}+\tau
B_{3}\\
T_{2}B_{2}-T_{1}B_{1}=0.\\
\end{array}\right.
\end{equation}
Noting that $T\times B=-N$ we also have
\begin{equation}\label{00E37}
T_{2}B_{1}-T_{1}B_{2}=N_{3}=0.
\end{equation}
Thus, we have
\begin{equation}\label{00E38}
\left\{\begin{array}{rl} T_{2}N_{2}-T_{1}N_{1}=-kT_{3}+\tau
B_{3},\;\;\;\langle1\rangle\\
 T_{2}^{2}-T_{1}^{2}+T_{3}'=0,\;\;\;\langle2\rangle \\
T_{2}B_{2}-T_{1}B_{1}=0,\;\;\langle3\rangle\\
 T_{2}B_{1}-T_{1}B_{2}=0.\;\;\langle4\rangle\\
\end{array}\right.
\end{equation}
Case A: $T_1^2\neq T_2^2$. In this case, Equations $\langle
3\rangle$ and $\langle 4\rangle$ in System (\ref{00E38}) viewed as
equations in $B_1$ and $B_2$ has a unique solution $B_1=B_2=0$. This
implies that $T_3=\langle N\times B, e_3\rangle=0$. Substitute this
into $\langle 2\rangle$ of System (\ref{00E38}) we have $T_1^2=
T_2^2$, a contradiction. Thus, we must have\\
Case B: $T_1^2= T_2^2$. In this case, Equation $\langle 2\rangle$ of
System (\ref{00E38}) implies that $T_3= {\rm constant}$. To
understand the meaning of this, we represent the unit tangent vector
$T$ as $T=\sin\alpha\cos\beta\, e_{1}+\sin\alpha\sin\beta
\,e_{2}+\cos\alpha \,e_{3} $, where $\alpha=\alpha(s),\;
\beta=\beta(s)$. With this representation, $T_3= {\rm constant}$
implies that $\cos\alpha={\rm constant}$ and hence $\alpha
(s)=\alpha_{0},\;{\rm a \; \;constant}$. This, together with $T_1^2=
T_2^2$, gives
\begin{equation}\label{OPM}
\sin\alpha_{0}(\cos\beta\pm\sin\beta)=0.
\end{equation}
If $\sin\alpha_{0}=0$, then we have $T_{1}=T_{2}=0$, and it follows
from the first equation of (\ref{00E34}) that $\nabla_{T}T=0$ which
means that $\gamma$ is a geodesic, a contradiction. Thus, we must
have $\sin\alpha_{0}\ne 0$, which, together with (\ref{OPM}),
implies that
\begin{equation}\notag
\cos\beta=\pm\sin\beta=\pm \frac{\sqrt{2}}{2},
\end{equation}
and hence,
\begin{equation}\label{00E529}T_{1}=\pm T_{2}=\pm
\frac{\sqrt{2}}{2}\sin\alpha_{0}.
\end{equation}
We use the first equation of (\ref{00E34}) again to get
\begin{equation}\notag
\nabla_{T}T =\sin\alpha_{0} \cos\alpha_{0}(\pm
\frac{\sqrt{2}}{2}e_{1}\mp \frac{\sqrt{2}}{2} e_{2})=kN
\end{equation}
which yields
\begin{equation}\label{L538}
N_{1}=\pm\frac{\sqrt{2}}{2},N_{2}= \mp\frac{\sqrt{2}}{2}
\end{equation}
since $k=|\nabla_{T}T|=|\sin\alpha_{0} \cos\alpha_{0}|$. By the
assumption that $\gamma$ is non-geodesic, we may assume, without
loss of generality, that $\sin\alpha_{0} \cos\alpha_{0}
>0$, so
\begin{equation}\label{ABCQ}
k=\sin\alpha_{0} \cos\alpha_{0}.
\end{equation}
Using Equations (\ref{00E529}), (\ref{L538}) and the fact that
$B=T\times N$ we have
\begin{equation}\label{L541}
\begin{array}{rl}
B_1= {\rm constant},\\
B_2= {\rm constant},\\
B_{3}=T_{1}N_{2}-T_{2}N_{1}=\mp\sin\alpha_{0}.
\end{array}
\end{equation}
It follows from (\ref{L541}), $\langle 3\rangle$ of (\ref{00E38}),
and the third equation of (\ref{00E34}) that
\begin{equation}\label{ABCQQ}
\tau^2=|\nabla_{T}B|^2= (T_1^2+T_2^2)B_3^2=\sin^4\alpha_{0}.
\end{equation}

Substituting (\ref{ABCQ}), (\ref{L541}) and (\ref{ABCQQ}) into the
third equation in (\ref{00E32}) we have
\begin{equation}\notag
\sin^2\alpha_{0}\cos^2
\alpha_{0}+\sin^4\alpha_{0}=2\sin^2\alpha_{0}-1,
\end{equation}
which implies
\begin{equation}\notag
\sin^2\alpha_{0} =1,
\end{equation}
and hence
\begin{equation}\notag
\cos^2\alpha_{0} =0.
\end{equation}
It follows that $\cos\alpha_{0} =0$ from which and (\ref{ABCQ}) we
conclude that $k=0$, i.e., $\gamma$ is a geodesic, a
contradiction.\\

Combining the the results in Cases A and B we complete the proof of
the theorem.
\end{proof}

\section{Linear biharmonic maps into Sol and Nil spaces}

In this section, we first derive the biharmonic map equation in
local coordinates and then we use it to classify linear biharmonic
maps from Euclidean space into Sol and Nil spaces. In the rest of
the paper, we adopt the Einstein summation convention that a
repeated upper and lower index means the summation of that index
over its range that is understood from the context.\\
{\bf 3.1 Biharmonic map equation in local coordinates}. \\
\begin{lemma}
Let $\varphi :(M^{m}, g)\longrightarrow (N^{n}, h)$ with $\varphi
(x^{1},\ldots, x^{m})=(\varphi^{1}(x), \ldots, \varphi^{n}(x))$ be a
map between Riemannian manifolds. With respect to local coordinates
$(x^{i})$ in $M$ and $(y^{\alpha})$ in $N$, $\varphi$ is biharmonic
if and only if it is a solution of the following system of PDE's
\begin{eqnarray}\label{BI2}
&&   g^{ij}(\tau^{\sigma}_{ij}
+\tau^{\alpha}_{j}\varphi^{\beta}_{i}{\bar
\Gamma_{\alpha\beta}^{\sigma}} +\frac{\partial}{\partial
x^{i}}(\tau^{\alpha}\varphi ^{\beta}_{j}{\bar
\Gamma_{\alpha\beta}^{\sigma}}) +\tau^{\alpha}\varphi
^{\beta}_{j}\varphi^{\rho}_{i}{\bar \Gamma_{\alpha\beta}^{\nu}}{\bar
\Gamma_{\nu\rho}^{\sigma}}\\\notag &&
-\Gamma^{k}_{ij}(\tau^{\sigma}_{k}+\tau^{\alpha}\varphi^{\beta}_{k}{\bar
\Gamma_{\alpha\beta}^{\sigma}}) -\tau^{\nu}{\varphi^{\alpha}}_{i}
\varphi^{\beta}_{j}{\bar R}_{\beta\,\alpha
\nu}^{\sigma})=0,\;\;\sigma=1,\, 2,\, \ldots, n.
\end{eqnarray}
\end{lemma}
\begin{proof}
Let $g_{ij}$ and $\Gamma_{ij}^{k}$ (resp. $h_{\alpha \beta}$ and
$\bar{\Gamma}^{\sigma}_{\alpha \beta}$) denote the components of
metric and the connection coefficients in the domain (resp. the
target) manifold with respect to the chosen local coordinates and
the natural frame $\{\frac{\partial}{\partial x^{i}}\}$ (resp.
$\{\frac{\partial}{\partial y^{\alpha}}\}$). The bitension field of
$\varphi$ can be computed as

\begin{eqnarray}\notag
\tau^{2}(\varphi)&=&  g^{ij} \left(
\nabla^{\varphi}_{\frac{\partial}{\partial
x^{i}}}\nabla^{\varphi}_{\frac{\partial}{\partial x^{j}}} -
\nabla^{\varphi}_{\nabla_{\frac{\partial}{\partial
x^{i}}}\frac{\partial}{\partial x^{j}}}\right)
\left(\tau(\varphi)\right) - {\rm Trace}_{g}{{\bar
R}(d\varphi,\tau(\varphi))d\varphi} \\\notag &= & g^{ij}
\nabla^{\varphi}_{\frac{\partial}{\partial
x^{i}}}\nabla^{\varphi}_{\frac{\partial}{\partial x^{j}}} \left(
\tau(\varphi)\right) -  g^{ij} \Gamma^{k}_{ij}
\nabla^{\varphi}_{\frac{\partial}{\partial x^{k}}}
 \left( \tau(\varphi)\right) \\\label{TAU}
&& \hskip 4.5cm -  g^{ij}{\bar R}( d\varphi
(\frac{\partial}{\partial x^{i}}), \tau(\varphi))d\varphi
(\frac{\partial}{\partial x^{j}})\\\notag &=& g^{ij}\left(
\nabla^{\varphi}_{\frac{\partial}{\partial
x^{i}}}\nabla^{\varphi}_{\frac{\partial}{\partial x^{j}}} \left(
\tau(\varphi)\right) -  \Gamma^{k}_{ij}
\nabla^{\varphi}_{\frac{\partial}{\partial x^{k}}}
 \left( \tau(\varphi)\right) - {\varphi^{\alpha}}_{i} \varphi^{\beta}_{j}
{\bar R}(\frac{\partial}{\partial y^{\alpha}},
\tau(\varphi))\frac{\partial}{\partial y^{\beta}}\right).
\end{eqnarray}
A direct computation gives
\begin{eqnarray}\label{T1}
\nabla^{\varphi}_{\frac{\partial}{\partial x^{k}}}
 \left( \tau(\varphi)\right)=\nabla^{\varphi}_{\frac{\partial}{\partial x^{k}}}
 (\tau^{\nu}\frac{\partial}{\partial
y^{\nu}})=\left(
\tau^{\sigma}_{k}+\tau^{\alpha}\varphi^{\beta}_{k}{\bar
\Gamma_{\alpha\beta}^{\sigma}} \right)\frac{\partial}{\partial
y^{\sigma}},
\end{eqnarray}
\begin{eqnarray}\notag
\nabla^{\varphi}_{\frac{\partial}{\partial
x^{i}}}\nabla^{\varphi}_{\frac{\partial}{\partial x^{j}}} \left(
\tau(\varphi)\right)&=&\nabla^{\varphi}_{\frac{\partial}{\partial
x^{i}}}\nabla^{\varphi}_{\frac{\partial}{\partial x^{j}}} (
\tau^{\nu}\frac{\partial}{\partial
y^{\nu}})=\nabla^{\varphi}_{\frac{\partial}{\partial x^{i}}}\left(
\tau^{\alpha}_{j}+\tau^{\nu}\varphi^{\beta}_{j}{\bar
\Gamma_{\beta\nu}^{\alpha}} \right)\frac{\partial}{\partial
y^{\alpha}}\\\label{T3} &=& \left( \tau^{\sigma}_{ij}
+\tau^{\alpha}_{j}\varphi^{\beta}_{i}{\bar
\Gamma_{\alpha\beta}^{\sigma}} +\frac{\partial}{\partial
x^{i}}(\tau^{\alpha}\varphi ^{\beta}_{j}{\bar
\Gamma_{\alpha\beta}^{\sigma}}) +\tau^{\alpha}\varphi
^{\beta}_{j}\varphi^{\rho}_{i}{\bar \Gamma_{\alpha\beta}^{\nu}}{\bar
\Gamma_{\nu\rho}^{\sigma}} \right)\frac{\partial}{\partial
y^{\sigma}},
\end{eqnarray}
and \begin{eqnarray}\label{T2} {\bar R}(\frac{\partial}{\partial
y^{\alpha}}, \tau(\varphi))\frac{\partial}{\partial y^{\beta}}=
\tau^{\nu}{\bar R}_{\beta\,\alpha
\nu}^{\sigma}\frac{\partial}{\partial y^{\sigma}}.
\end{eqnarray}
 Substitute Equations (\ref{T1}), (\ref{T3}), and (\ref{T2}) into
(\ref{TAU}) we have

\begin{eqnarray}\notag
\tau^{2}(\varphi) &=&   g^{ij}(\tau^{\sigma}_{ij}
+\tau^{\alpha}_{j}\varphi^{\beta}_{i}{\bar
\Gamma_{\alpha\beta}^{\sigma}} +\frac{\partial}{\partial
x^{i}}(\tau^{\alpha}\varphi ^{\beta}_{j}{\bar
\Gamma_{\alpha\beta}^{\sigma}}) +\tau^{\alpha}\varphi
^{\beta}_{j}\varphi^{\rho}_{i}{\bar \Gamma_{\alpha\beta}^{\nu}}{\bar
\Gamma_{\nu\rho}^{\sigma}}\\
&&
-\Gamma^{k}_{ij}(\tau^{\sigma}_{k}+\tau^{\alpha}\varphi^{\beta}_{k}{\bar
\Gamma_{\alpha\beta}^{\sigma}}) -\tau^{\nu}{\varphi^{\alpha}}_{i}
\varphi^{\beta}_{j}{\bar R}_{\beta\,\alpha
\nu}^{\sigma})\frac{\partial}{\partial y^{\sigma}},
\end{eqnarray}
from which the lemma follows.
\end{proof}

When the domain manifold is a Euclidean space then we have
\begin{corollary}\label{COR1}
Let $\varphi :\mathbb{R}^{m}\longrightarrow (N^{n}, h)$ with
$\varphi (x^{1},\ldots, x^{m})=(\varphi^{1}(x), \ldots,
\varphi^{n}(x))$ be a map from a Euclidean space into a Riemannian
manifold. Then $\varphi$ is biharmonic if and only if it is a
solution of the following system of PDE's
\begin{eqnarray}\label{BihE}
&&  \Delta \tau^{\sigma} +\langle \nabla \tau^{\alpha}, \nabla
\varphi^{\beta}\rangle{\bar \Gamma_{\alpha\beta}^{\sigma}} +\langle
\nabla \varphi ^{\beta},\nabla (\tau^{\alpha}{\bar
\Gamma_{\alpha\beta}^{\sigma}})\rangle\\\notag && +\, \langle \nabla
\varphi ^{\beta}, \nabla \varphi^{\rho}\rangle \tau^{\alpha}{\bar
\Gamma_{\alpha\beta}^{\nu}}{\bar \Gamma_{\nu\rho}^{\sigma}}
 -\tau^{\nu}\langle \nabla \varphi^{\alpha}, \nabla
\varphi^{\beta}\rangle{\bar R}_{\beta\,\alpha
\nu}^{\sigma}=0,\;\;\sigma=1,\, 2,\, \ldots, n.
\end{eqnarray}
\end{corollary}

{\bf  Linear biharmonic maps into Sol space}. Let
$(\mathbb{R}^{3},g_{Sol})$ denote Sol space, where the metric can be
written as $g_{Sol}=e^{2y^{3}}({\rm d}y^{1})^{2}+e^{-2y^{3}}({\rm
d}y^{2})^{2}+({\rm d}y^{3})^{2}$ with respect to the standard
coordinates $(y^{1},y^{2},y^{3})$ in $\mathbb{R}^{3}$. Then a direct
computation gives the following components of Sol metric and the
coefficients of the connection:
\begin{align*}\nonumber
& g_{11}=e^{2y^{3}},\; g_{22}=e^{-2y^{3}},\;g_{33}=1,\; {\rm
all\;\;other\;\;} g_{ij}=0;\\\nonumber
&g^{11}=e^{-2y^{3}},\;g^{22}=e^{2y^{3}},\;g^{33}=1,\; {\rm
all\;\;other\;\;} g^{ij}=0;\\
\end{align*}
\begin{equation}\label{Q758}
\begin{array}{lll}
\bar{\Gamma}^{1}_{11}=\bar{\Gamma}^{2}_{11}=0,\bar{\Gamma}^{3}_{11}=-e^{2y^{3}};\\
\bar{\Gamma}^{1}_{12}=0, \bar{\Gamma}^{2}_{12}=0,\bar{\Gamma}^{3}_{12}=0;\\
\bar{\Gamma}^{1}_{13}=1, \bar{\Gamma}^{2}_{13}=0,\bar{\Gamma}^{3}_{13}=0;\\
\bar{\Gamma}^{1}_{21}=0, \bar{\Gamma}^{2}_{21}=0,\bar{\Gamma}^{3}_{21}=0;\\
\bar{\Gamma}^{1}_{22}=\bar{\Gamma}^{2}_{22}=0,\bar{\Gamma}^{3}_{22}=e^{-2y^{3}};\\
\bar{\Gamma}^{1}_{23}=0, \bar{\Gamma}^{2}_{23}=-1,\bar{\Gamma}^{3}_{23}=0;\\
\bar{\Gamma}^{1}_{31}=1, \bar{\Gamma}^{2}_{31}=0,\bar{\Gamma}^{3}_{31}=0;\\
\bar{\Gamma}^{1}_{32}=0, \bar{\Gamma}^{2}_{32}=-1,\bar{\Gamma}^{3}_{32}=0;\\
\bar{\Gamma}^{1}_{33}=0, \bar{\Gamma}^{2}_{33}=0,\bar{\Gamma}^{3}_{33}=0.\\
\end{array}
\end{equation}
By our convention of curvature operator and the following notation
for the components of the Riemannian curvature
\begin{equation}
\bar{R}(\frac{\partial}{\partial y^{i}},\frac{\partial}{\partial
y^{j}})\frac{\partial}{\partial y^{k}}
=\bar{R}_{k\,ij}^{l}\frac{\partial}{\partial y^{l}}
\end{equation}
we have
\begin{equation}\label{Riem}
\bar{R}_{k\,ij}^{l}=\frac{\partial}{\partial
y^{i}}\bar{\Gamma}_{kj}^{l}-\frac{\partial}{\partial
y^{j}}\bar{\Gamma}_{ki}^{l}+\bar{\Gamma}_{it}^{l}\bar{\Gamma}_{kj}^{t}-\bar{\Gamma}_{jt}^{l}
\bar{\Gamma}_{ki}^{t}
\end{equation}
A straightforward computation using (\ref{Q758}) and (\ref{Riem})
gives the following components of the Riemannian curvature of Sol
space:
\begin{equation}\label{Q790}
\left\{\begin{array}{rl} {\bar R}_{221}^{1}=-e^{-2y^{3}},{\bar
R}_{331}^{1}=1,
{\bar R}_{212}^{1}=e^{-2y^{3}},{\bar R}_{313}^{1}=-1;\\
{\bar R}_{121}^{2}=e^{2y^{3}},{\bar R}_{112}^{2}=-e^{2y^{3}},
{\bar R}_{332}^{2}=1,{\bar R}_{323}^{2}=-1;\\
{\bar R}_{131}^{3}=-e^{2y^{3}},{\bar R}_{232}^{3}=-e^{-2y^{3}},
{\bar R}_{113}^{3}=e^{2y^{3}},{\bar R}_{223}^{3}=e^{-2y^{3}}.\\
\end{array}\right.
\end{equation}

Now we are ready to prove the following classification theorem for
linear biharmonic maps into Sol space.\\
\begin{theorem}\label{MAI1}
Let $\varphi:\mathbb{R}^{m}\longrightarrow(\mathbb{R}^{3}, g_{Sol})$
with
\begin{equation*}\label{MAP51}
\varphi (x)=\left(\begin{array}{ccccc} a^{1}_{1} & a^{1}_{2} &
\cdots & a^{1}_{m}
\\
a^{2}_{1} & a^{2}_{2} & \cdots
& a^{2}_{m}\\
a^{3}_{1} & a^{3}_{2} & \cdots & a^{3}_{m}
\end{array}\right)
\left(\begin{array}{ccccc}
x^{1}\\

x^{2}\\
\vdots\\
x^{m}

\end{array}\right),
\end{equation*}
i.e., $\varphi (X)=(A^{1}X^{t}, A^{2}X^{t}, A^{3}X^{t})$ be a linear
map into Sol space, where $A^{i}$ denotes the row vectors of the
representation matrix. Then, $\varphi$ is a biharmonic map if and
only if it is a harmonic map, which is equivalent to either $ (i)\;
A^{3}=0,\; |A^{1}|^{2}=|A^{2}|^{2}$, or, $(ii)\; A^{3}\ne 0,\;
A^{1}=A^{2}=0$.
\end{theorem}
\begin{proof}
With respect to the standard Cartesian coordinates $(x^{i})$ in
$\mathbb{R}^{m}$ and $(y^{\alpha})$ in $\mathbb{R}^{3}$, the tension
field of $\varphi$ is given by
\begin{equation}\label{Q651}
\begin{array}{lll}
\tau(\varphi):={\rm Trace}_{g}(\nabla {\rm d}\varphi)\in
\Gamma(\varphi^{-1}TN)\\  =\tau^{\sigma}\frac{\partial}{\partial
y^{\sigma}}\\
=g^{ij}(\varphi_{ij}^{\sigma}-\Gamma_{ij}^{k}\varphi_{k}^{\sigma}+
\bar{\Gamma}^{\sigma}_{\alpha
\beta}\varphi^{\alpha}_{i}\varphi^{\beta}_{j})\frac{\partial}{\partial
y^{\sigma}}\\=(\sum\limits_{i=1}^{m}\bar{\Gamma}^{\sigma}_{\alpha
\beta}\varphi^{\alpha}_{i}\varphi^{\beta}_{i})\frac{\partial}{\partial
y^{\sigma}}\\=\bar{\Gamma}^{\sigma}_{\alpha \beta}A^{\alpha}\cdot
A^{\beta}\frac{\partial}{\partial y^{\sigma}},
\end{array}
\end{equation}
where and in the sequel, $A^{\alpha}\cdot A^{\beta}$ denotes the
inner product and $|A^{\alpha}|$ the norm of the vectors in Euclidean space.\\
Putting $\tau(\varphi)=\tau^{\sigma}\frac{\partial}{\partial
y^{\sigma}}$ and substituting (\ref{Q758}) into Equation
(\ref{Q651}) we find the following components of the tension field
of $\varphi$

\begin{equation}\label{Q896}
\begin{array}{lll}
\tau^{1}=2A^{1}\cdot A^{3},\\
\tau^{2}=-2A^{2}\cdot A^{3},\\
\tau^{3}=|A^{2}|^{2}e^{-2y^{3}}-|A^{1}|^{2}e^{2y^{3}},
\end{array}
\end{equation}
where and in the sequel $y^{3}=A^3 X^t$.\\

A further computation gives

\begin{eqnarray}\notag
&& \left( \Delta \tau^{\sigma} +\langle \nabla \tau^{\alpha}, \nabla
\varphi^{\beta}\rangle{\bar
\Gamma_{\alpha\beta}^{\sigma}}\right)\frac{\partial}{\partial
y^{\sigma}}=\\\notag &&\left( \Delta \tau^{\sigma} +A^{\beta}\cdot
\nabla \tau^{\alpha}{\bar
\Gamma_{\alpha\beta}^{\sigma}}\right)\frac{\partial}{\partial
y^{\sigma}}=\\\label{Q904} &&-2A^{1}\cdot
A^{3}(|A^{2}|^{2}e^{-2y^{3}}+|A^{1}|^{2}e^{2y^{3}}
)\frac{\partial}{\partial y^{1}}\\\notag &&+2A^{2}\cdot
A^{3}(|A^{2}|^{2}e^{-2y^{3}}+|A^{1}|^{2}e^{2y^{3}}
)\frac{\partial}{\partial y^{2}}\\\notag &&+4|
A^{3}|^{2}(|A^{2}|^{2}e^{-2y^{3}}-|A^{1}|^{2}e^{2y^{3}}
)\frac{\partial}{\partial y^{3}},
\end{eqnarray}

\begin{eqnarray}\notag
&&\langle \nabla \varphi ^{\beta},\nabla (\tau^{\alpha}{\bar
\Gamma_{\alpha\beta}^{\sigma}})\rangle\frac{\partial}{\partial
y^{\sigma}}\\\notag &&=( A ^{\beta}\cdot \nabla \tau^{\alpha}{\bar
\Gamma_{\alpha\beta}^{\sigma}}+ \tau^{\alpha}A^{\beta}\cdot \nabla
{\bar \Gamma_{\alpha\beta}^{\sigma}} )\frac{\partial}{\partial
y^{\sigma}}\\\label{BBCA} &&=-2A^{1}\cdot
A^{3}(|A^{2}|^{2}e^{-2y^{3}}
+|A^{1}|^{2}e^{2y^{3}})\frac{\partial}{\partial y^{1}}\\\notag
&&+2A^{2}\cdot A^{3}(|A^{2}|^{2}e^{-2y^{3}}
+|A^{1}|^{2}e^{2y^{3}})\frac{\partial}{\partial y^{2}}\\\notag
&&+[-4(A^{1}\cdot A^{3})^{2}e^{2y^{3}} +4(A^{2}\cdot
A^{3})^{2}e^{-2y^{3}}]\frac{\partial}{\partial y^{3}},
\end{eqnarray}

\begin{equation}\label{W2441}
\begin{array}{lll}
\langle \nabla \varphi ^{\beta}, \nabla \varphi^{\rho}\rangle
\tau^{\alpha}{\bar \Gamma_{\alpha\beta}^{\nu}}{\bar
\Gamma_{\nu\rho}^{\sigma}}\frac{\partial}{\partial
y^{\sigma}}=\tau^{\alpha} A ^{\beta}\cdot A^{\rho} {\bar
\Gamma_{\alpha\beta}^{\nu}}{\bar
\Gamma_{\nu\rho}^{\sigma}}\frac{\partial}{\partial
y^{\sigma}}\\
=[2A^{1}\cdot A^{3}|A^{3}|^{2} +A^{1}\cdot
A^{3}|A^{2}|^{2}e^{-2y^{3}}-3A^{1}\cdot A^{3}|A^{1}|^{2}e^{2y^{3}}
\\-2A^{2}\cdot A^{3}A^{1}\cdot A^{2}e^{-2y^{3}}]\frac{\partial}{\partial y^{1}}
+[-2A^{2}\cdot A^{3}|A^{3}|^{2} -A^{2}\cdot
A^{3}|A^{1}|^{2}e^{2y^{3}}\\+3A^{2}\cdot A^{3}|A^{2}|^{2}e^{-2y^{3}}
+2A^{1}\cdot A^{3}A^{1}\cdot
A^{2}e^{2y^{3}}]\frac{\partial}{\partial y^{2}} +[-2(A^{1}\cdot
A^{3})^{2}e^{2y^{3}} \\+|A^{1}|^{4}e^{4y^{3}}+2(A^{2}\cdot
A^{3})^{2}e^{-2y^{3}}
-|A^{2}|^{4}e^{-4y^{3}}]\frac{\partial}{\partial y^{3}},
\end{array}
\end{equation}
 and
\begin{equation}\label{B3031}
\begin{array}{lll}
\tau^{\nu}\langle \nabla \varphi^{\alpha}, \nabla
\varphi^{\beta}\rangle{\bar R}_{\beta\,\alpha
\nu}^{\sigma}\frac{\partial}{\partial y^{\sigma}}=A^{\beta}\cdot
A^{\alpha}\tau^{\nu}{\bar R}_{\beta\,\alpha
\nu}^{\sigma}\frac{\partial}{\partial
y^{\sigma}}\\
=[-3A^{1}\cdot A^{3}|A^{2}|^{2}e^{-2y^{3}} -2A^{2}\cdot
A^{3}A^{1}\cdot A^{2}e^{-2y^{3}}\\ +A^{1}\cdot
A^{3}|A^{1}|^{2}e^{2y^{3}}+2A^{1}\cdot
A^{3}|A^{3}|^{2}]\frac{\partial}{\partial y^{1}} +[2A^{1}\cdot
A^{3}A^{1}\cdot A^{2}e^{2y^{3}} \\+3A^{2}\cdot
A^{3}|A^{1}|^{2}e^{2y^{3}} -A^{2}\cdot
A^{3}|A^{2}|^{2}e^{-2y^{3}}-2A^{2}\cdot
A^{3}|A^{3}|^{2}]\frac{\partial}{\partial y^{2}}\\
+[-|A^{1}|^{4}e^{4y^{3}}+|A^{2}|^{4}e^{-4y^{3}} -2(A^{1}\cdot
A^{3})^{2}e^{2y^{3}} +2(A^{2}\cdot
A^{3})^{2}e^{-2y^{3}})]\frac{\partial}{\partial y^{3}}.
\end{array}
\end{equation}
It follows from Equations (\ref{Q904}), (\ref{BBCA}), (\ref{W2441}),
(\ref{B3031}) and Corollary \ref{COR1} that the linear map $\varphi$
is a biharmonic map if and only if
\begin{equation}\label{B368}
\left\{\begin{array}{rl}-8A^{1}\cdot
A^{3}|A^{1}|^{2}e^{2y^{3}}=0\\
8A^{2}\cdot A^{3}|A^{2}|^{2}e^{-2y^{3}}=0\\
4(|A^{2}|^{2}|A^{3}|^{2}+(A^{2}\cdot A^{3})^{2})e^{-2y^{3}}-
4(|A^{1}|^{2}|A^{3}|^{2}+(A^{1}\cdot A^{3})^{2})e^{2y^{3}}\\
+2|A^{1}|^{4}e^{4y^{3}} -2|A^{2}|^{4}e^{-4y^{3}}=0.
\end{array}\right.
\end{equation}
Solving System of equations (\ref{B368}) we have either $ (i)\;
A^{3}=0,\; |A^{1}|^{2}=|A^{2}|^{2}$, or $(ii)\; A^{3}\ne 0,\;
A^{1}=A^{2}=0$. It follows from Equation (\ref{Q896}) that in both
cases the tension field vanishes identically, i.e., $\varphi$ is
also harmonic. Therefore, we obtain the theorem.
\end{proof}

{\bf Linear biharmonic maps into Nil space}. Let
$(\mathbb{R}^{3},g_{Nil})$ denote Nil space, where the metric with
respect to the standard coordinates $(y^{1},y^{2},y^{3})$ in
$\mathbb{R}^{3}$ can be written as $g_{Nil}=({\rm d}y^{1})^{2}+({\rm
d}y^{2})^{2}+({\rm d}y^{3}-y^{1}{\rm d}y^{2})^{2}$. Then an easy
computation gives the following components of Nil metric and the
coefficients of its connection:
\begin{align}\nonumber
& g_{11}=1,\; g_{12}=g_{13}=0,\;
g_{22}=1+(y^{1})^{2},\;g_{23}=-y^{1},\;g_{33}=1;\\\nonumber
&g^{11}=1,\;g^{12}=g^{13}=0,\;g^{22}=1,\;g^{23}=y^{1},\;g^{33}=1+(y^{1})^{2};\\\notag
\end{align}
\begin{equation}\label{Q812}
\begin{array}{lll}
\bar{\Gamma}^{1}_{11}=\bar{\Gamma}^{2}_{11}=0,\bar{\Gamma}^{3}_{11}=0;\\
\bar{\Gamma}^{1}_{12}=0, \bar{\Gamma}^{2}_{12}=\frac{y^{1}}{2},\bar{\Gamma}^{3}_{12}=\frac{(y^{1})^{2}-1}{2};\\
\bar{\Gamma}^{1}_{13}=0, \bar{\Gamma}^{2}_{13}=-\frac{1}{2},\bar{\Gamma}^{3}_{13}=-\frac{y^{1}}{2};\\
\bar{\Gamma}^{1}_{21}=0, \bar{\Gamma}^{2}_{21}=\frac{y^{1}}{2},\bar{\Gamma}^{3}_{21}=\frac{(y^{1})^{2}-1}{2};\\
\bar{\Gamma}^{1}_{22}=-y^{1},\bar{\Gamma}^{2}_{22}=0,\bar{\Gamma}^{3}_{22}=0;\\
\bar{\Gamma}^{1}_{23}=\frac{1}{2}, \bar{\Gamma}^{2}_{23}=0,\bar{\Gamma}^{3}_{23}=0;\\
\bar{\Gamma}^{1}_{31}=0, \bar{\Gamma}^{2}_{31}=-\frac{1}{2},\bar{\Gamma}^{3}_{31}=-\frac{y^{1}}{2};\\
\bar{\Gamma}^{1}_{32}=\frac{1}{2}, \bar{\Gamma}^{2}_{32}=0,\bar{\Gamma}^{3}_{32}=0;\\
\bar{\Gamma}^{1}_{33}=0, \bar{\Gamma}^{2}_{33}=0,\bar{\Gamma}^{3}_{33}=0.\\
\end{array}
\end{equation}

A further computation using (\ref{Riem}) and (\ref{Q812}) gives the
following components of the Riemannian curvature of Nil space:
\begin{equation}\label{Q845}
\left\{\begin{array}{rl} {\bar
R}_{212}^{1}=-\frac{3}{4}+\frac{(y^{1})^{2}}{4},{\bar
R}_{213}^{1}=-\frac{y^{1}}{4},{\bar
R}_{221}^{1}=\frac{3}{4}-\frac{(y^{1})^{2}}{4},{\bar R}_{312}^{1}=
-\frac{y^{1}}{4},\\{\bar R}_{231}^{1}=\frac{y^{1}}{4}, {\bar
R}_{313}^{1}=\frac{1}{4},{\bar
R}_{321}^{1}=\frac{y^{1}}{4},{\bar R}_{331}^{1}=-\frac{1}{4};\\
{\bar R}_{112}^{2}=\frac{3}{4},{\bar R}_{121}^{2}=-\frac{3}{4},
{\bar R}_{223}^{2}=-\frac{y^{1}}{4},\\{\bar
R}_{232}^{2}=\frac{y^{1}}{4},{\bar R}_{323}^{2}=\frac{1}{4},{\bar
R}_{332}^{2}=-\frac{1}{4}
;\\
{\bar R}_{112}^{3}=y^{{1}},{\bar R}_{113}^{3}=-\frac{1}{4}, {\bar
R}_{121}^{3}=-y^{{1}},{\bar R}_{131}^{3}=\frac{1}{4},\\ {\bar
R}_{223}^{3}=-\frac{(y^{1})^{2}+1}{4},{\bar
R}_{232}^{3}=\frac{(y^{1})^{2}+1}{4},
{\bar R}_{323}^{3}=\frac{y^{1}}{4},{\bar R}_{332}^{3}=-\frac{y^{1}}{4}.\\
\end{array}\right.
\end{equation}

\begin{theorem}\label{MAI2}
Let $\varphi:\mathbb{R}^{m}\longrightarrow(\mathbb{R}^{3}, g_{Nil})$
with
\begin{equation*}\label{MAP51}
\varphi (x)=\left(\begin{array}{ccccc} a^{1}_{1} & a^{1}_{2} &
\cdots & a^{1}_{m}
\\
a^{2}_{1} & a^{2}_{2} & \cdots
& a^{2}_{m}\\
a^{3}_{1} & a^{3}_{2} & \cdots & a^{3}_{m}
\end{array}\right)
\left(\begin{array}{ccccc}
x^{1}\\

x^{2}\\
\vdots\\
x^{m}

\end{array}\right),
\end{equation*}
i.e., $\varphi (X)=(A^{1}X^{t}, A^{2}X^{t}, A^{3}X^{t})$ be a linear
map into Nil space, where $A^{i}$ denotes the row vectors of the
representation matrix. Then, $\varphi$ is a biharmonic map if and
only if it is a harmonic map, which is equivalent to either $ (i)\;
A^{1}=0,\; A^{2}\cdot A^{3}=0$, or, $ (ii)\; A^{1}=0,\; A^{2}=0$,
or, $(iii)\; A^{2}= 0,\; A^{1} \cdot A^{3}=0$.
\end{theorem}

\begin{proof}
Taking the standard Cartesian coordinates $(x^{i})$ in
$\mathbb{R}^{m}$, $(y^{\alpha})$ in $\mathbb{R}^{3}$ and
substituting (\ref{Q812}) into (\ref{Q651}) we find the tension
field of $\varphi$ to be

\begin{equation}\label{M1072}
\begin{array}{lll}
\tau(\varphi)=\tau^{\sigma}\frac{\partial}{\partial
y^{\sigma}}\\
=[-|A^{2}|^{2}y^{1}+A^{2}\cdot A^{3}]\frac{\partial}{\partial
y^{1}}\\
+[A^{1}\cdot A^{2}y^{1}-A^{1}\cdot A^{3}]\frac{\partial}{\partial
y^{2}}\\+[A^{1}\cdot A^{2}(y^{1})^{2} -A^{1}\cdot
A^{3}y^{1}-A^{1}\cdot A^{2}]\frac{\partial}{\partial
y^{3}},\\
\end{array}
\end{equation}
or,
\begin{equation}\label{M1168}
\begin{array}{lll}
\tau^{1}=-|A^{2}|^{2}y^{1}+A^{2}\cdot A^{3},\\
\tau^{2}=A^{1}\cdot A^{2}y^{1}-A^{1}\cdot A^{3},\\
\tau^{3}=A^{1}\cdot A^{2}(y^{1})^{2} -A^{1}\cdot
A^{3}y^{1}-A^{1}\cdot A^{2},
\end{array}
\end{equation}
where and in the sequel $y^{1}=A^1 X^t$.\\

A straightforward computation yields
\begin{equation}\label{M168}
\begin{array}{lll}
\Delta \tau^{\sigma}\frac{\partial}{\partial y^{\sigma}}
=2A^{1}\cdot A^{2}|A^{1}|^{2}\frac{\partial}{\partial
y^{3}},\\
\end{array}
\end{equation}

\begin{equation}\label{M239}
\begin{array}{lll}
\langle \nabla \tau^{\alpha}, \nabla \varphi^{\beta}\rangle{\bar
\Gamma_{\alpha\beta}^{\sigma}}\frac{\partial}{\partial
y^{\sigma}}=A^{\beta}\cdot \nabla \tau^{\alpha}{\bar
\Gamma_{\alpha\beta}^{\sigma}}\frac{\partial}{\partial
y^{\sigma}}\\
= [-\frac{1}{2}A^{1}\cdot A^{2}(|A^{1}|^{2}+|A^{2}|^{2})y^{{1}}
+\frac{1}{2}A^{1}\cdot A^{3}(|A^{1}|^{2}+|A^{2}|^{2})]\frac{\partial}{\partial y^{2}}\\
+[-\frac{1}{2}A^{1}\cdot A^{2}(|A^{1}|^{2}+|A^{2}|^{2})(y^{{1}})^{2}
+\frac{1}{2}A^{1}\cdot A^{3}(|A^{1}|^{2}+|A^{2}|^{2})y^{{1}}\\
+\frac{1}{2}A^{1}\cdot
A^{2}(|A^{2}|^{2}-|A^{1}|^{2})]\frac{\partial}{\partial y^{3}},
\end{array}
\end{equation}

\begin{equation}\label{M315}
\begin{array}{lll}
\langle \nabla \varphi ^{\beta},\nabla (\tau^{\alpha}{\bar
\Gamma_{\alpha\beta}^{\sigma}})\rangle\frac{\partial}{\partial
y^{\sigma}}\\=( A ^{\beta}\cdot \nabla \tau^{\alpha}{\bar
\Gamma_{\alpha\beta}^{\sigma}}+ \tau^{\alpha}A^{\beta}\cdot \nabla
{\bar \Gamma_{\alpha\beta}^{\sigma}} )\frac{\partial}{\partial
y^{\sigma}}\\ =[-(A^{1}\cdot A^{2})^{2}y^{{1}} +(A^{1}\cdot
A^{2})(A^{1}\cdot A^{3})]\frac{\partial}{\partial
y^{1}}\\+[-A^{1}\cdot A^{2}|A^{2}|^{2}y^{{1}}
+\frac{1}{2}(A^{1}\cdot A^{2})(A^{2}\cdot A^{3})+
\frac{1}{2}A^{1}\cdot A^{3}|A^{2}|^{2}]\frac{\partial}{\partial y^{2}}\\
+[-\frac{3}{2}A^{1}\cdot A^{2}|A^{2}|^{2}(y^{{1}})^{2} +(A^{1}\cdot
A^{2})(A^{2}\cdot A^{3})y^{{1}} \\+A^{1}\cdot
A^{3}|A^{2}|^{2}y^{{1}} +\frac{1}{2}A^{1}\cdot
A^{2}|A^{2}|^{2}-\frac{1}{2}(A^{1}\cdot A^{3})(A^{2}\cdot
A^{3})]\frac{\partial}{\partial y^{3}},
\end{array}
\end{equation}
\begin{equation}\label{M458}
\begin{array}{lll}
\langle \nabla \varphi ^{\beta}, \nabla \varphi^{\rho}\rangle
\tau^{\alpha}{\bar \Gamma_{\alpha\beta}^{\nu}}{\bar
\Gamma_{\nu\rho}^{\sigma}}\frac{\partial}{\partial
y^{\sigma}}\\\tau^{\alpha} A ^{\beta}\cdot A^{\rho} {\bar
\Gamma_{\alpha\beta}^{\nu}}{\bar
\Gamma_{\nu\rho}^{\sigma}}\frac{\partial}{\partial
y^{\sigma}}\\=[\frac{1}{4}|A^{2}|^{4}(y^{{1}})^{3}
-\frac{3}{4}A^{2}\cdot A^{3}|A^{2}|^{2}(y^{{1}})^{2}\\
+\frac{1}{2}(A^{2}\cdot A^{3})^{2}y^{{1}} -\frac{1}{2}(A^{1}\cdot
A^{2})^{2}y^{{1}}+\frac{1}{4}|A^{2}|^{2} |A^{3}|^{2}y^{{1}}
+\frac{1}{4}|A^{2}|^{4}y^{{1}}\\+\frac{1}{2}(A^{1}\cdot
A^{2})(A^{1}\cdot A^{3}) -\frac{1}{4}A^{2}\cdot A^{3}|
A^{3}|^{2}-\frac{1}{4}A^{2}\cdot
A^{3}|A^{2}|^{2}]\frac{\partial}{\partial
y^{1}}\\+[-\frac{1}{4}A^{1}\cdot A^{2}|A^{2}|^{2}(y^{{1}})^{3}
+\frac{1}{4}A^{1}\cdot A^{3}|A^{2}|^{2}(y^{{1}})^{2}
\\+\frac{1}{2}(A^{2}\cdot A^{3})(A^{1}\cdot A^{2})(y^{{1}})^{2}
-\frac{1}{2}A^{1}\cdot A^{2}|A^{2}|^{2}y^{{1}}\\
+\frac{1}{4}A^{1}\cdot A^{2}|A^{1}|^{2}y^{{1}}-\frac{1}{4}A^{1}\cdot
A^{2}|A^{3}|^{2}y^{{1}}-\frac{1}{2}(A^{1}\cdot A^{3})(A^{2}\cdot
A^{3})y^{{1}}\\
+\frac{1}{2}(A^{1}\cdot A^{2})(A^{2}\cdot A^{3})
+\frac{1}{4}A^{1}\cdot A^{3}|A^{3}|^{2}-\frac{1}{4}A^{1}\cdot
A^{3}|A^{1}|^{2}]\frac{\partial}{\partial y^{2}}\\
+[-\frac{1}{4}A^{1}\cdot A^{2}|A^{2}|^{2}(y^{{1}})^{4}
+\frac{1}{4}A^{1}\cdot
A^{3}|A^{2}|^{2}(y^{{1}})^{3}+\frac{1}{2}(A^{2}\cdot
A^{3})(A^{1}\cdot A^{2})(y^{{1}})^{3}\\
 +\frac{1}{4}A^{1}\cdot
A^{2}|A^{1}|^{2}(y^{{1}})^{2} -\frac{1}{2}(A^{2}\cdot
A^{3})(A^{1}\cdot A^{3})(y^{{1}})^{2} -\frac{1}{4}A^{1}\cdot
A^{2}|A^{3}|^{2}(y^{{1}})^{2}\\-\frac{1}{2}A^{1}\cdot
A^{3}|A^{2}|^{2}y^{{1}} +\frac{1}{4}A^{1}\cdot
A^{3}|A^{3}|^{2}y^{{1}}-\frac{1}{4}A^{1}\cdot
A^{3}|A^{1}|^{2}y^{{1}}\\
+\frac{1}{2}(A^{2}\cdot A^{3})(A^{1}\cdot A^{3})
+\frac{1}{4}A^{1}\cdot A^{2}|A^{2}|^{2} -\frac{1}{4}A^{1}\cdot
A^{2}|A^{1}|^{2}]\frac{\partial}{\partial y^{3}},
\end{array}
\end{equation}
and
\begin{equation}\label{M335}
\begin{array}{lll}
\tau^{\nu}\langle \nabla \varphi^{\alpha}, \nabla
\varphi^{\beta}\rangle{\bar R}_{\beta\,\alpha
\nu}^{\sigma}\frac{\partial}{\partial y^{\sigma}}=A^{\beta}\cdot
A^{\alpha}\tau^{\nu}{\bar R}_{\beta\,\alpha
\nu}^{\sigma}\frac{\partial}{\partial y^{\sigma}}\\=[\frac{1}{4}|
A^{2}|^{4}(y^{{1}})^{3} -\frac{3}{4}A^{2}\cdot
A^{3}|A^{2}|^{2}(y^{{1}})^{2} -\frac{1}{2}(A^{1}\cdot
A^{2})^{2}y^{{1}}+\frac{1}{2}(A^{2}\cdot A^{3})^{2}y^{{1}}
\\-\frac{3}{4}|A^{2}|^{4}y^{{1}}+\frac{1}{4}|A^{2}|^{2}|A^{3}|^{2}y^{{1}}
+\frac{1}{2}(A^{1}\cdot A^{2})(A^{1}\cdot
A^{3})+\frac{3}{4}A^{2}\cdot
A^{3}|A^{2}|^{2}\\-\frac{1}{4}A^{2}\cdot
A^{3}|A^{3}|^{2}]\frac{\partial}{\partial
y^{1}}+[-\frac{1}{4}A^{1}\cdot A^{2}|A^{2}|^{2}(y^{{1}})^{3}
+\frac{1}{4}A^{1}\cdot A^{3}|A^{2}|^{2}(y^{{1}})^{2}
\\+\frac{1}{2}(A^{2}\cdot A^{3})(A^{1}\cdot A^{2})(y^{{1}})^{2}
+\frac{3}{4}A^{1}\cdot A^{2}|A^{1}|^{2}y^{{1}}+A^{1}\cdot
A^{2}|A^{2}|^{2}y^{{1}} \\-\frac{1}{2}(A^{2}\cdot A^{3})(A^{1}\cdot
A^{3})y^{{1}} -\frac{1}{4}A^{1}\cdot A^{2}|A^{3}|^{2}y^{{1}}
-\frac{3}{4}A^{1}\cdot A^{3}|A^{1}|^{2} \\-(A^{1}\cdot
A^{2})(A^{2}\cdot A^{3})+\frac{1}{4}A^{1}\cdot
A^{3}|A^{3}|^{2}]\frac{\partial}{\partial y^{2}}
+[-\frac{1}{4}A^{1}\cdot A^{2}|A^{2}|^{2}(y^{{1}})^{4}
\\+\frac{1}{4}A^{1}\cdot
A^{3}|A^{2}|^{2}(y^{{1}})^{3}+\frac{1}{2}(A^{2}\cdot
A^{3})(A^{1}\cdot A^{2})(y^{{1}})^{3}
 +\frac{3}{4}A^{1}\cdot
A^{2}|A^{1}|^{2}(y^{{1}})^{2} \\+A^{1}\cdot
A^{2}|A^{2}|^{2}(y^{{1}})^{2} -\frac{1}{2}(A^{2}\cdot
A^{3})(A^{1}\cdot A^{3})(y^{{1}})^{2}-\frac{1}{4}A^{1}\cdot
A^{2}|A^{3}|^{2}(y^{{1}})^{2} \\-\frac{3}{4}A^{1}\cdot
A^{3}|A^{1}|^{2}y^{{1}}- (A^{1}\cdot A^{2})(A^{2}\cdot A^{3})y^{{1}}
+\frac{1}{4}A^{1}\cdot A^{3}|A^{3}|^{2}y^{{1}}\\
+\frac{1}{4}A^{1}\cdot A^{2}|A^{1}|^{2}+\frac{1}{4}A^{1}\cdot
A^{2}|A^{2}|^{2}]\frac{\partial}{\partial y^{3}}.
\end{array}
\end{equation}

By Equations (\ref{M168}), (\ref{M239}), (\ref{M315}), (\ref{M458}),
(\ref{M335}) and Corollary \ref{COR1} we conclude that $\varphi$ is
biharmonic if and only if

\begin{equation}\label{M665}
\begin{array}{lll}-(A^{1}\cdot A^{2})^{2}y^{{1}}
+|A^{2}|^{4}y^{{1}}\\+(A^{1}\cdot A^{2})(A^{1}\cdot A^{3})
-A^{2}\cdot A^{3}| A^{2}|^{2}
=0,\\
\end{array}
\end{equation}
\begin{equation}\label{M683}
\begin{array}{lll}-A^{1}\cdot
A^{2}|A^{1}|^{2}y^{{1}}-3A^{1}\cdot A^{2}|A^{2}|^{2}y^{{1}}
\\+A^{1}\cdot A^{3}|A^{2}|^{2}+A^{1}\cdot A^{3}|A^{1}|^{2}\\
+2(A^{1}\cdot A^{2})(A^{2}\cdot A^{3})=0,\\
\end{array}
\end{equation}
and
\begin{equation}\label{M703}
\begin{array}{lll}
-A^{1}\cdot A^{2}|A^{1}|^{2}(y^{{1}})^{2} -3A^{1}\cdot
A^{2}|A^{2}|^{2}(y^{{1}})^{2} +A^{1}\cdot
A^{3}|A^{1}|^{2}y^{{1}}\\+A^{1}\cdot A^{3}|A^{2}|^{2}y^{{1}}
+2(A^{1}\cdot A^{2})(A^{2}\cdot A^{3})y^{{1}}\\+A^{1}\cdot
A^{2}|A^{1}|^{2}+ A^{1}\cdot A^{2}|A^{2}|^{2}
=0.\\
\end{array}
\end{equation}
To solve the System of biharmonic map equations (\ref{M665}),
(\ref{M683}), and (\ref{M703}) we consider the following two cases.\\
Case I: $A^{1}=0$. Noting that $y^{1}=A^1 X^t=0$ we use Equation
(\ref{M665}) to have\\
\begin{equation}\notag
-A^{2}\cdot A^{3}|A^{2}|^{2}=0,
\end{equation}
which implies either $A^{2}\cdot A^{3}=0$, or $ |A^{2}|^{2}=0
\;\;i.e. \;\;A^{2}=0$.

Substituting $A^{1}=0, \; A^{2}\cdot A^{3}=0$ or $A^{1}=0, \;
A^{2}=0$ into (\ref{M683}) and (\ref{M703}) we see that they are
both solutions of biharmonic map equations.\\

Case II: $A^{1}\neq 0$. Note that in this case, $y^{1}=A^1 X^t \neq
0$. We can view system of biharmonic map equations (\ref{M665}),
(\ref{M683}), and (\ref{M703})  as a system of polynomial equations
in $y^1$. It is not difficult to check that
$A^2=0,\; A^1\cdot A^3=0$ is the only solution in this case.\\

Combining Case I and II we obtain the last statement of the theorem.
A direct checking using the tension field equation (\ref{M1072}) we
see that all these maps are also harmonic maps. Thus, we complete
the proof of the theorem.
\end{proof}

\end{document}